\documentclass[oneside,english]{amsart}
\usepackage[T1]{fontenc}
\usepackage[latin9]{inputenc}
\usepackage{amstext}
\usepackage{amsthm}
\usepackage{amssymb}
\usepackage{mathdots}
\usepackage{esint}

\makeatletter
\numberwithin{equation}{section}
\numberwithin{figure}{section}
\theoremstyle{plain}
\newtheorem{thm}{\protect\theoremname}
\theoremstyle{remark}
\newtheorem{rem}[thm]{\protect\remarkname}
\theoremstyle{definition}
\newtheorem{defn}[thm]{\protect\definitionname}
\theoremstyle{plain}
\newtheorem{prop}[thm]{\protect\propositionname}

\usepackage{breqn}

\makeatother

\usepackage{babel}
\providecommand{\definitionname}{Definition}
\providecommand{\propositionname}{Proposition}
\providecommand{\remarkname}{Remark}
\providecommand{\theoremname}{Theorem}

\begin{document}
\title{Symplectic aspects of the tt{*}-Toda equations}
\author{Ryosuke Odoi}
\address{Department of Pure and Applied Mathematics, Faculty of Science and
Engineering, Waseda University, 3-4-1 Okubo, Shinjuku, Tokyo 169-8555
JAPAN}
\email{ryosuke.odoi@moegi.waseda.jp}
\begin{abstract}
We evaluate explicitly, in terms of the asymptotic data, the ratio
of the constant pre-factors in the large and small $x$ asymptotics
of the tau functions for global solutions of the tt{*}-Toda equations.
This constant problem for the sinh-Gordon equation, which is the
case $n=1$ of the tt{*}-Toda equations, was solved by C.\ A.\ Tracy
\cite{Tr}. We also introduce natural symplectic structures on the
space of asymptotic data and on the space of monodromy data for a
wider class of solutions, and show that these symplectic structures
are preserved by the Riemann-Hilbert correspondence.
\end{abstract}

\maketitle
\tableofcontents{}

\section{Introduction}

Painlev\'e equations may be formulated as Hamiltonian systems. This
has led to an important role in the theory of such equations for concepts
from classical mechanics and symplectic geometry, such as canonical
coordinates, tau-functions, and moduli spaces of solutions with symplectic
structures. The benefit of the symplectic point of view is that it
illuminates a path to the study of more general nonlinear differential
equations, especially those which are ``integrable''.

The Painlev\'e equations themselves are scalar ordinary differential
equations of second order, and this facilitates explicit calculations.
For systems, or higher-order equations, geometry plays a more essential
role.

Indeed, a considerable amount of general theory has been developed,
for example by Hitchin \cite{H}, Boalch \cite{B}, building on earlier
work of Schlesinger \cite{Sch}, Jimbo-Miwa-Ueno \cite{JMU}. On
the other hand examples are rather scarce, partly because of the difficulty
of carrying out explicit calculations, and partly because of the lack
of interesting concrete example for higher rank systems.

The purpose of this article is to explain some symplectic aspects
of the tt{*}-Toda equations, a system of nonlinear ordinary differential
equation of ``Painlev\'e type'', which is a relatively recent example
arising in physics. The tt{*} equations (topological--anti topological
fusion equations) arose in the work of Cecotti and Vafa on supersymmetric
quantum field theory (\cite{CV}, \cite{D}), and the tt{*}-Toda equations
are a special case of these equations of ``Toda type''. The simplest
nontrivial case of the tt{*}-Toda equations is the (radial) sinh-Gordon
equation, which is in fact a case of the third Painlev\'e equation.
It was investigated -- for similar physical reasons -- by McCoy-Tracy-Wu
\cite{MTW}, and their work had far-reaching consequences.

More recently, the tt{*}-Toda equations were investigated in detail
by Guest-Its-Lin \cite{II,III} and by Mochizuki \cite{M1,M2}, and
our motivation was to put some of these results into a symplectic
context and investigate them further. We have succeeded to do this
only for a certain subset of solutions -- an open subset of the moduli
space of all solutions -- but the results are encouraging, and have
already led to a new application, which we shall explain later.

The paper is organized as follows. After a brief review of the tt{*}-Toda
equations in section 2, we give their Hamiltonian formulation in section
3. In section 4 we explain the symplectic structures on the space
of solutions that we consider and on a corresponding space of monodromy
data. The correspondence between solutions (asymptotic data) and monodromy
data (Stokes matrices and connection matrices) is an example of the
Riemann-Hilbert correspondence for meromorphic connections with irregular
singularities. Our first main result (Theorem \ref{thm:symp}) is
that this correspondence preserves the symplectic structures. This
is consistent with the general results of Boalch \cite{B}, but we
shall go further and give an explicit generating function which relates
the corresponding canonical coordinates (Theorem \ref{thm:gen_fcn}).

In section 5 we give an application of these results to the asymptotics
of the tau function. For each solution of the tt{*}-Toda equation
there is a corresponding tau function, and it is the properties these
tau functions (rather than the solutions themselves) which are important
for many applications in physics.

A recent example is the work of Its, Lisovyy, and Tykhyy \cite{ILT},
in which the structure of the tau functions was elucidated using representation
theory (conformal blocks), as a consequence of AGT duality in physics.
This was used to solve the ``constant problem'' for Painlev\'e
equations, i.e. the problem of finding the constant which relates
the short-distance and long-distance expansions of the tau function.
In the case of the (radial) sine-Gordon equation, a rigorous and more
direct proof was given by Its and Prokhorov \cite{IP}. Our second
main result makes use of their method in order to solve this ``constant
problem'' for the tt{*}-Toda equations. As in \cite{IP}, we find
that the (explicit) generating function plays a crucial role (Theorem
\ref{thm:const}).

This work is part of the author's Ph.\ D.\ thesis at Waseda University.
He would like to acknowledge his supervisor, Prof.\ Martin Guest
for his support throughout this work and his guidance in this field.
He would like to acknowledge Prof.\ Alexander Its for his friendly
advice and his suggestions regarding the constant problem. He is also
glad to acknowledge the financial support from the Mathematics and
Physics Unit ``Multiscale Analysis, Modelling and Simulation'',
Top Global University Project, Waseda University.

\section{The tt{*}-Toda equations}

Let a positive integer $n$ be fixed. The tt{*}-Toda equations are
\begin{equation}
2\left(w_{i}\right)_{t\bar{t}}=-e^{2(w_{i+1}-w_{i})}+e^{2(w_{i}-w_{i-1})},\:w_{i}:\mathbb{C}^{*}\rightarrow\mathbb{R},\:i\in\mathbb{Z},\label{eq:tt}
\end{equation}
where, for all i, $w_{i}=w_{i+n+1}$, $w_{i}=w_{i}(\left|t\right|)$
$(t\in\mathbb{C}^{*})$, and
\[
w_{0}+w_{n}=0,w_{1}+w_{n-1}=0,\;\dots\;\text{(anti-symmetry condition)}.
\]

The equations (\ref{eq:tt}) are equivalent to the flatness of $\nabla:=d+\alpha$,
i.e. the zero curvature equation $d\alpha+\alpha\wedge\alpha=0$,
where 
\[
\alpha:=\left(w_{t}+\frac{1}{\lambda}W^{T}\right)dt+\left(-w_{\bar{t}}+\lambda W\right)d\bar{t},
\]
\[
w=\begin{pmatrix}w_{0}\\
 & \ddots\\
 &  & w_{n}
\end{pmatrix},\,W=\begin{pmatrix}0 & e^{w_{1}-w_{0}}\\
 & 0 & \ddots\\
 &  & \ddots & e^{w_{n}-w_{n-1}}\\
e^{w_{0}-w_{n}} &  &  & 0
\end{pmatrix}.
\]

The tt{*}-Toda equations are also equivalent to the isomonodromy condition
for the following ordinary differential equation
\begin{equation}
\frac{d\Psi}{d\zeta}=\left(-\frac{1}{\zeta^{2}}W-\frac{1}{\zeta}xw_{x}+x^{2}W^{T}\right)\Psi,\label{eq:aux}
\end{equation}
 where $x:=\left|t\right|$.

Generically, the local solutions near $x=0$ of the tt{*}-Toda equations
are parametrized by real numbers $\gamma_{i},\rho_{i}$ as follows
\cite{III,M1,M2}:
\begin{equation}
2w_{i}(x)=\gamma_{i}\log x+\rho_{i}+o(1)\quad\text{as }x\rightarrow0.\label{eq:o(1)}
\end{equation}
 We call the parameters $\gamma_{i},\rho_{i}$ the asymptotic data.
``Generically'' means $-2<\gamma_{i+1}-\gamma_{i}<2$; the general
case has $-2\leq\gamma_{i+1}-\gamma_{i}\leq2$. We assume the generic
condition from now on.

There is another important set of data $m_{i}$, $e_{i}^{\mathbb{R}}$
called the monodromy data. These are eigenvalues of certain matrices
$M$ and $E$, which are related to monodromy data such as Stokes
matrices. See \cite{III} or the appendix for details. The proof in
\cite{III} is for the case $n=3$, but exactly the same method provides
the results of Theorem \ref{thm:GIL} and \ref{thm:global-solutions}
below for general $n$.
\begin{thm}
\label{thm:GIL}\cite{III}The monodromy data $m_{i}$, $e_{i}^{\mathbb{R}}$
may be expressed in terms of the asymptotic data as follows:
\begin{align*}
m_{i} & =-\frac{1}{2}\gamma_{i}\\
e_{i}^{\mathbb{R}} & =\begin{cases}
e^{\rho_{i}}2^{2\gamma_{i}}\frac{X_{n-i}(\gamma_{0},\dots,\gamma_{(n-1)/2},-\gamma_{(n-1)/2},\dots,-\gamma_{0})}{X_{i}(\gamma_{0},\dots,\gamma_{(n-1)/2},-\gamma_{(n-1)/2},\dots,-\gamma_{0})} & n:\text{odd}\\
e^{\rho_{i}}2^{2\gamma_{i}}\frac{X_{n-i}(\gamma_{0},\dots,\gamma_{(n-2)/2},0,-\gamma_{(n-2)/2},\dots,-\gamma_{0})}{X_{i}(\gamma_{0},\dots,\gamma_{(n-2)/2},0,-\gamma_{(n-2)/2},\dots,-\gamma_{0})} & n:\text{even}
\end{cases}
\end{align*}
where 
\[
X_{k}(\gamma_{0},\dots,\gamma_{n}):=\prod_{j=1}^{n}\Gamma(\frac{\gamma_{k}-\gamma_{k+j}+2j}{2(n+1)})\;(\gamma_{j+n+1}=\gamma_{j}).
\]
\end{thm}

Global solutions can be parametrized only by the $\gamma_{i}$ (or
only by the $m_{i}$), that is, for global solutions the $\rho_{i}$
are determined by the $\gamma_{i}$:
\begin{thm}
\label{thm:global-solutions}\cite{III}For global solutions (i.e.
solutions which are smooth for $0<x<\infty$) we have
\[
\rho_{i}=-(2\log2)\gamma_{i}+\log(X_{i}/X_{n-i}),
\]
i.e. $e_{i}^{\mathbb{R}}=1$.
\end{thm}

\section{The Hamiltonian formulation\label{sec:The-Hamiltonian-formulation}}

Next, we introduce a Hamiltonian function and a symplectic form.

Let $\left\lfloor x\right\rfloor :=\max\{n\in\mathbb{Z}:n\leq x\}$
for $x\in\mathbb{R}$. The tt{*}-Toda equations can be written as
a non-autonomous Hamiltonian system, 
\begin{align}
\left(w_{i}\right)_{x} & =\frac{\partial H}{\partial\tilde{w}_{i}}=\frac{\tilde{w}_{i}}{x}\label{eq:*}\\
\left(\tilde{w}_{i}\right)_{x} & =-\frac{\partial H}{\partial w_{i}}=-2x\left(e^{2(w_{i+1}-w_{i})}-e^{2(w_{i}-w_{i-1})}\right),\label{eq:**}
\end{align}
 on the phase space $\mathbb{R}^{2\left\lfloor (n-1)/2\right\rfloor +2}=\{(w,\tilde{w})\}$
($w=(w_{0},\dots,w_{\left\lfloor (n-1)/2\right\rfloor })$, $\tilde{w}=(\tilde{w}_{0},\dots,\tilde{w}_{\left\lfloor (n-1)/2\right\rfloor })$)
equipped with the symplectic structure 
\[
\theta:=\sum_{i=0}^{\left\lfloor (n-1)/2\right\rfloor }dw_{i}\wedge d\tilde{w}_{i}
\]
 where the Hamiltonian $H$ is defined by 
\[
H(w,\tilde{w};x):=\frac{1}{2x}\sum_{i=0}^{\left\lfloor (n-1)/2\right\rfloor }\tilde{w}_{i}^{2}-x\sum_{i=1}^{\left\lfloor (n-1)/2\right\rfloor }e^{2(w_{i}-w_{i-1})}-\frac{x}{2}\left(e^{-4w_{\left\lfloor (n-1)/2\right\rfloor }}+e^{4w_{0}}\right).
\]

The symplectic form $\theta$ is asymptotic to $\sum_{i=0}^{\left\lfloor (n-1)/2\right\rfloor }d(\rho_{i}/2)\wedge d(\gamma_{i}/2)$
as $x\rightarrow0$.
\begin{rem}
The Hamiltonian system may be written in terms of $X:=\log x$ as
follows:
\[
H(w,\tilde{w};X):=\frac{1}{2e^{X}}\sum_{i=0}^{\left\lfloor (n-1)/2\right\rfloor }\tilde{w}_{i}^{2}-e^{X}\sum_{i=1}^{\left\lfloor (n-1)/2\right\rfloor }e^{2(w_{i}-w_{i-1})}-\frac{e^{X}}{2}\left(e^{-4w_{\left\lfloor (n-1)/2\right\rfloor }}+e^{4w_{0}}\right).
\]
\begin{align*}
\left(w_{i}\right)_{X} & =\frac{\partial e^{X}H}{\partial\tilde{w}_{i}}=\tilde{w}_{i}\\
\left(\tilde{w}_{i}\right)_{X} & =-\frac{\partial e^{X}H}{\partial w_{i}}=-2e^{2X}\left(e^{2(w_{i+1}-w_{i})}-e^{2(w_{i}-w_{i-1})}\right).
\end{align*}
\end{rem}

\section{Symplectic structures and the Riemann-Hilbert correspondence}

Both the asymptotic data $\gamma_{i}$, $\rho_{i}$ and the monodromy
data $m_{i}$, $\log e_{i}^{\mathbb{R}}$ can be considered as defining
local charts of the moduli space of solutions. From Theorem \ref{thm:GIL}
we can show that the transformation between two charts via the Riemann-Hilbert
correspondence is symplectic with respect to the ``obvious'' symplectic
structure. The symplectic form $2\theta$ we define in section \ref{sec:The-Hamiltonian-formulation}
is asymptotic to the left hand side of the equality below as $x\rightarrow0$.
\begin{thm}
\label{thm:symp}
\[
-\frac{1}{2}\sum_{i=0}^{\left\lfloor (n-1)/2\right\rfloor }d\gamma_{i}\wedge d\rho_{i}=\sum_{i=0}^{\left\lfloor (n-1)/2\right\rfloor }dm_{i}\wedge d\log e_{i}^{\mathbb{R}}.
\]
\end{thm}

\begin{rem}
\label{rem:KKAH}The left hand side is related to the Kirillov-Kostant
form on a coadjoint orbit, and the right hand side is related to the
Atiyah-Hitchin form on the space of the based rational maps of degree
$n+1$ from $\mathbb{C}P^{1}$ to itself. Thus, both symplectic forms
arise naturally from geometry. We shall present details of these facts
elsewhere.
\end{rem}

Theorem \ref{thm:symp} can be verified by direct calculation, but
we prefer to give a proof by showing the existence of a generating
function. The generating function will play an important role later.
\begin{defn}
\label{def:generating function}Let 
\begin{align}
F(\rho_{0},\dots,\rho_{\left\lfloor (n-1)/2\right\rfloor },m_{0},\dots,m_{\left\lfloor (n-1)/2\right\rfloor }): & =-\sum_{i=0}^{\left\lfloor (n-1)/2\right\rfloor }\rho_{i}m_{i}+2\log2\sum_{i=0}^{\left\lfloor (n-1)/2\right\rfloor }m_{i}^{2}\label{eq:F}\\
 & +\frac{n+1}{2}\sum_{k=0}^{n}\sum_{j=1}^{n}\psi^{(-2)}\left(\frac{m_{k-j}-m_{k}+j}{n+1}\right)\nonumber 
\end{align}
 where $m_{j+n+1}=m_{j}$ and $m_{j}=-m_{n-j}$. Here $\psi^{(-2)}(z)=\int_{0}^{z}\log\Gamma(x)dx=\frac{z(1-z)}{2}+\frac{z}{2}\log2\pi+z\log\Gamma(z)-\log G(1+z)$,
and $G$ is the Barnes G-function.
\end{defn}

\begin{thm}
\label{thm:gen_fcn}The function $F$ is a generating function of
the transformation 
\[
(m_{0},\dots,m_{\left\lfloor (n-1)/2\right\rfloor },\rho_{0},\dots,\rho_{\left\lfloor (n-1)/2\right\rfloor })\mapsto(m_{0},\dots,m_{\left\lfloor (n-1)/2\right\rfloor },\log e_{0}^{\mathbb{R}},\dots,\log e_{\left\lfloor (n-1)/2\right\rfloor }^{\mathbb{R}})
\]
 with respect to the given symplectic forms. More precisely, $F$
satisfies 
\[
m_{i}=-\frac{\partial F}{\partial\rho_{i}},\;\log e_{i}^{\mathbb{R}}=-\frac{\partial F}{\partial m_{i}}.
\]
\end{thm}

\begin{proof}
The first identity is obvious. We show the second identity. Let 
\[
\tilde{K}(m_{0},\dots,m_{n}):=(n+1)\sum_{i=0}^{n}\sum_{j=1}^{n}\psi^{(-2)}(\frac{m_{i}-m_{i+j}+j}{n+1}),
\]
where $m_{j+n+1}=m_{j}$. Let 
\[
K(m_{0},\dots,m_{\left\lfloor (n-1)/2\right\rfloor }):=\frac{1}{2}\tilde{K}(m_{0},\dots,m_{\left\lfloor (n-1)/2\right\rfloor },-m_{\left\lfloor (n-1)/2\right\rfloor },\dots,-m_{0}).
\]
This $K$ is the last term of $F$ in (\ref{eq:F}). From the definition
of $\log e_{i}^{\mathbb{R}}$ and $F$, it suffices to show that 
\begin{align*}
 & \frac{\partial K}{\partial m_{k}}(m_{0},\dots,m_{\left\lfloor (n-1)/2\right\rfloor })\\
 & =\log\left(\frac{X_{k}(m_{0},\dots.m_{\left\lfloor (n-1)/2\right\rfloor },-m_{\left\lfloor (n-1)/2\right\rfloor },\dots,-m_{0})}{X_{n-k}(m_{0},\dots.m_{\left\lfloor (n-1)/2\right\rfloor },-m_{\left\lfloor (n-1)/2\right\rfloor },\dots,-m_{0})}\right).
\end{align*}
We can easily obtain that $X_{n-k}(-m_{n},\dots,-m_{0})=\prod_{j=1}^{n}\Gamma(\frac{-m_{k}+m_{k-j}+j}{n+1})$
and that $\tilde{K}(m_{0},\dots,m_{n})=(n+1)\sum_{k=0}^{n}\sum_{j=1}^{n}\psi^{(-2)}(\frac{m_{k-j}-m_{k}+j}{n+1})$.
Then we obtain 
\[
\frac{\partial\tilde{K}}{\partial m_{k}}=\log(X_{k}(m_{0},\dots,m_{n})/X_{n-k}(-m_{n},\dots,-m_{0})).
\]

Hence we have
\begin{align*}
 & \frac{\partial K}{\partial m_{k}}(m_{0},\dots,m_{\left\lfloor (n-1)/2\right\rfloor })\\
= & \frac{1}{2}\biggl(\frac{\partial\tilde{K}}{\partial m_{k}}(m_{0},\dots,m_{\left\lfloor (n-1)/2\right\rfloor },-m_{\left\lfloor (n-1)/2\right\rfloor },\dots,-m_{0})\\
 & -\frac{\partial\tilde{K}}{\partial m_{n-k}}(m_{0},\dots,m_{\left\lfloor (n-1)/2\right\rfloor },-m_{\left\lfloor (n-1)/2\right\rfloor },\dots,-m_{0})\biggr)\\
= & \frac{1}{2}\bigl(\log X_{k}(m_{0},\dots.m_{\left\lfloor (n-1)/2\right\rfloor },-m_{\left\lfloor (n-1)/2\right\rfloor },\dots,-m_{0})\\
 & -\log X_{n-k}(m_{0},\dots.m_{\left\lfloor (n-1)/2\right\rfloor },-m_{\left\lfloor (n-1)/2\right\rfloor },\dots,-m_{0})\\
 & -\log X_{n-k}(m_{0},\dots.m_{\left\lfloor (n-1)/2\right\rfloor },-m_{\left\lfloor (n-1)/2\right\rfloor },\dots,-m_{0})\\
 & +\log X_{k}(m_{0},\dots.m_{\left\lfloor (n-1)/2\right\rfloor },-m_{\left\lfloor (n-1)/2\right\rfloor },\dots,-m_{0})\bigr)\\
= & \log X_{k}(m_{0},\dots.m_{\left\lfloor (n-1)/2\right\rfloor },-m_{\left\lfloor (n-1)/2\right\rfloor },\dots,-m_{0})\\
 & -\log X_{n-k}(m_{0},\dots.m_{\left\lfloor (n-1)/2\right\rfloor },-m_{\left\lfloor (n-1)/2\right\rfloor },\dots,-m_{0}).
\end{align*}
 This completes the proof.
\end{proof}

\section{Tau functions and the constant problem}

In this section, we assume for simplicity that $n=3$, so $w=(w_{0},w_{1},w_{2},w_{3})$
with $w_{2}=-w_{1}$, $w_{3}=-w_{0}$. For general $n$ the same method
applies. We consider only the global solutions, i.e., we assume $\log e_{i}^{\mathbb{R}}=0$,
which means that the $\rho_{i}$'s are determined by the $\gamma_{j}$'s,
as in Theorem \ref{thm:global-solutions}. The following calculation
is motivated by Theorem 1 in \cite{IP}. See also \cite{P} for further
details.
\begin{defn}
Let us define the tau function of a global solution $w$ by 
\[
\log\tau^{w}(x_{1},x_{2})=\intop_{x_{1}}^{x_{2}}H(w_{i}(x),\tilde{w}_{i}(x),x)dx
\]
 where $H$ is the Hamiltonian function.
\end{defn}

\begin{rem}
Usually the tau function is defined (up to a multiplicative constant)
by $\log\tau^{w}(x)=\intop^{x}H(w_{i}(x),\tilde{w}_{i}(x),x)dx$.
In that notation we have $\tau^{w}(x_{1},x_{2})=\tau^{w}(x_{2})/\tau^{w}(x_{1})$.
\end{rem}

The Hamiltonian function is

\[
H(x,w_{0},w_{1},\tilde{w}_{0},\tilde{w}_{1})=\frac{1}{2x}(\tilde{w}_{0}^{2}+\tilde{w}_{1}^{2})-xe^{2(w_{1}-w_{0})}-\frac{x}{2}\left(e^{-4w_{1}}+e^{4w_{0}}\right)
\]
 and is quasihomogeneous, that is,

\[
H(w,\lambda\tilde{w};\lambda x)=\lambda H(w,\tilde{w};x)\text{ for any }\lambda>0.
\]
 It follows that 
\[
\sum_{i=0}^{1}\tilde{w}_{i}\frac{\partial H}{\partial\tilde{w}_{i}}+x\frac{\partial H}{\partial x}=H.
\]

For the solution $(w_{0}(x),w_{1}(x),\tilde{w}_{0}(x),\tilde{w}_{1}(x))$
of (\ref{eq:*}) and (\ref{eq:**}), we have 
\begin{align*}
\sum_{i=0}^{1}\tilde{w}_{i}(x)\frac{\partial H}{\partial\tilde{w}_{i}}(x,w_{0}(x),w_{1}(x),\tilde{w}_{0}(x),\tilde{w}_{1}(x)) & =\sum_{i=0}^{1}\tilde{w}_{i}(x)\left(w_{i}\right)_{x}(x)\\
x\frac{\partial H}{\partial x}(x,w_{0}(x),w_{1}(x),\tilde{w}_{0}(x),\tilde{w}_{1}(x)) & =-H(w_{0}(x),w_{1}(x),\tilde{w}_{0}(x),\tilde{w}_{1}(x))+\\
 & \frac{dxH(x,w_{0}(x),w_{1}(x),\tilde{w}_{0}(x),\tilde{w}_{1}(x))}{dx}.
\end{align*}

The first equality is obvious. The second equality follows from $\frac{dxH}{dx}=x\frac{dH}{dx}+H$
and 
\begin{align*}
 & \frac{dH(x,w_{0}(x),w_{1}(x),\tilde{w}_{0}(x),\tilde{w}_{1}(x))}{dx}\\
 & =\frac{\partial H}{\partial x}(x,w_{0}(x),w_{1}(x),\tilde{w}_{0}(x),\tilde{w}_{1}(x))+\left(w_{i}\right)_{x}(x)\frac{\partial H}{\partial w_{i}}(x,w_{0}(x),w_{1}(x),\tilde{w}_{0}(x),\tilde{w}_{1}(x))\\
 & +\left(\tilde{w}_{i}\right)_{x}(x)\frac{\partial H}{\partial\tilde{w}_{i}}(x,w_{0}(x),w_{1}(x),\tilde{w}_{0}(x),\tilde{w}_{1}(x))\\
 & =\frac{\partial H}{\partial x}(x,w_{0}(x),w_{1}(x),\tilde{w}_{0}(x),\tilde{w}_{1}(x))\\
 & -\left(w_{i}\right)_{x}(x)\left(\tilde{w}_{i}\right)_{x}(x)+\left(\tilde{w}_{i}\right)_{x}(x)\left(w_{i}\right)_{x}(x)\\
 & =\frac{\partial H}{\partial x}(x,w_{0}(x),w_{1}(x),\tilde{w}_{0}(x),\tilde{w}_{1}(x)).
\end{align*}
Then it follows that 
\begin{prop}
$H=\tilde{w}_{0}\left(w_{0}\right)_{x}+\tilde{w}_{1}\left(w_{1}\right)_{x}-H+\frac{d}{dx}(xH).$
\end{prop}

Let $S(x_{1},x_{2}):=\intop_{x_{1}}^{x_{2}}\left(\sum_{i=0}^{1}\tilde{w}_{i}\left(w_{i}\right)_{x}-H\right)dx$,
which is called the classical action, the functional from which we
can derive the Euler-Lagrange equation using the fundamental lemma
of calculus of variations. We obtain 
\begin{align*}
\frac{\partial S(x_{1},x_{2})}{\partial\gamma_{j}} & =\intop_{x_{1}}^{x_{2}}\left(\sum_{i=0}^{1}\left(\left(\tilde{w}_{i}\right)_{\gamma_{j}}\left(w_{i}\right)_{x}+\tilde{w}_{i}\left(\left(w_{i}\right)_{x}\right)_{\gamma_{j}}\right)-\left(H\right)_{\gamma_{j}}\right)dx\\
 & =\intop_{x_{1}}^{x_{2}}\left(\sum_{i=0}^{1}\left(\left(\tilde{w}_{i}\right)_{\gamma_{j}}\frac{\partial H}{\partial\tilde{w}_{i}}+\tilde{w}_{i}\left(\left(w_{i}\right)_{x}\right)_{\gamma_{j}}\right)-\sum_{i=0}^{1}\left(\frac{\partial H}{\partial w_{i}}\left(w_{i}\right)_{\gamma_{j}}+\frac{\partial H}{\partial\tilde{w}_{i}}\left(\tilde{w}_{i}\right)_{\gamma_{j}}\right)\right)dx\\
 & =\intop_{x_{1}}^{x_{2}}\left(\sum_{i=0}^{1}\left(\left(\tilde{w}_{i}\right)_{\gamma_{j}}\frac{\partial H}{\partial\tilde{w}_{i}}-\left(\tilde{w}_{i}\right)_{x}\left(w_{i}\right)_{\gamma_{j}}\right)-\sum_{i=0}^{1}\left(\frac{\partial H}{\partial w_{i}}\left(w_{i}\right)_{\gamma_{j}}+\frac{\partial H}{\partial\tilde{w}_{i}}\left(\tilde{w}_{i}\right)_{\gamma_{j}}\right)\right)dx\\
 & +\left.\left(\sum_{i=0}^{1}\tilde{w}_{i}\left(w_{i}\right)_{\gamma_{j}}\right)\right|_{x_{1}}^{x_{2}}\\
 & =\left.\left(\tilde{w}_{0}\left(w_{0}\right)_{\gamma_{j}}+\tilde{w}_{1}\left(w_{1}\right)_{\gamma_{j}}\right)\right|_{x_{1}}^{x_{2}}.
\end{align*}

The second equality follows from (\ref{eq:*}) and the chain rule,
the third from integration by parts, and the fourth from (\ref{eq:**}).

From the proposition and the definition of the $\tau$ function, we
obtain 
\begin{align}
\frac{\partial}{\partial\gamma_{j}}\log\tau^{w}(x_{1},x_{2}) & =\frac{\partial}{\partial\gamma_{j}}\int_{x_{1}}^{x_{2}}\left(\sum_{i=0}^{1}\tilde{w}_{i}\left(w_{i}\right)_{x}-H+\frac{d}{dx}(xH)\right)dx\nonumber \\
 & =\frac{\partial S(x_{1},x_{2})}{\partial\gamma_{j}}+\left(x_{2}H(x_{2})-x_{1}H(x_{1})\right)_{\gamma_{j}}\nonumber \\
 & =\left.\left(\sum_{i=0}^{1}\tilde{w}_{i}\left(w_{i}\right)_{\gamma_{j}}\right)\right|_{x_{1}}^{x_{2}}+\left(x_{2}H(x_{2})-x_{1}H(x_{1})\right)_{\gamma_{j}}.\label{eq:gamma_tau}
\end{align}

At $x=0$ the form of (\ref{eq:o(1)}) is 
\begin{align*}
w_{i}(x) & =\frac{\gamma_{i}}{2}\log x+\frac{\rho_{i}}{2}+O(x^{\varepsilon_{i}}),\quad x\rightarrow0
\end{align*}
 for some $\varepsilon_{i}>0$ (which depends on $\gamma_{0}$ and
$\gamma_{1}$); this can be shown as in Theorem 14.1 of \cite{FIKN}
for the case $n=1$. This formula is differentiable in $x$ and the
$\gamma_{i}$'s. Therefore 
\begin{align*}
\tilde{w}_{i} & =\frac{\gamma_{i}}{2}+O(x^{\varepsilon_{i}}),\\
\left(w_{i}\right)_{\gamma_{j}} & =\frac{\delta_{i,j}}{2}\log x+\frac{1}{2}\left(\rho_{i}\right)_{\gamma_{j}}+O(x^{\varepsilon_{i}}\log x),\\
\left(\tilde{w}_{i}\right)_{\gamma_{j}} & =\frac{\delta_{i,j}}{2}+O(x^{\varepsilon_{i}}\log x)
\end{align*}
 as $x\rightarrow0$.

At $x=\infty$, from \cite{II}, if $s_{1}^{\mathbb{R}}\neq0$,
\begin{equation}
w_{i}(x)=-s_{1}^{\mathbb{R}}2^{-\frac{7}{4}}(\pi x)^{-\frac{1}{2}}e^{-2\sqrt{2}x}+O(x^{-1}e^{-2\sqrt{2}x})\quad\text{as }x\rightarrow\infty,\label{eq:infinity}
\end{equation}

where $s_{1}^{\mathbb{R}}=-2\cos\frac{\pi}{4}(\gamma_{0}+1)-2\cos\frac{\pi}{4}(\gamma_{1}+3)$.
If $s_{1}^{\mathbb{R}}=0$, we have
\begin{align*}
w_{0}(x) & =s_{2}^{\mathbb{R}}2^{-\frac{5}{2}}(\pi x)^{-\frac{1}{2}}e^{-4x}+O(x^{-1}e^{-4x})\sim O(x^{-1}e^{-2\sqrt{2}x})\\
w_{1}(x) & =-s_{2}^{\mathbb{R}}2^{-\frac{5}{2}}(\pi x)^{-\frac{1}{2}}e^{-4x}+O(x^{-1}e^{-4x})\sim O(x^{-1}e^{-2\sqrt{2}x}),
\end{align*}
 so the equation (\ref{eq:infinity}) holds for any generic $(\gamma_{0},\gamma_{1})$.
The equation (\ref{eq:infinity}) is also differentiable in $x$ and
the $\gamma_{i}$'s, so 
\begin{align*}
\tilde{w}_{i}(x) & =s_{1}^{\mathbb{R}}2^{-\frac{1}{4}}\sqrt{\pi}x^{\frac{1}{2}}e^{-2\sqrt{2}x}+O(e^{-2\sqrt{2}x}),\\
\left(w_{i}\right)_{\gamma_{j}} & =-\left(s_{1}^{\mathbb{R}}\right)_{\gamma_{j}}2^{-\frac{7}{4}}(\pi x)^{-\frac{1}{2}}e^{-2\sqrt{2}x}+O(x^{-1}e^{-2\sqrt{2}x}),\\
\left(\tilde{w}_{i}\right)_{\gamma_{j}} & =\left(s_{1}^{\mathbb{R}}\right)_{\gamma_{j}}2^{-\frac{1}{4}}\sqrt{\pi}x^{\frac{1}{2}}e^{-2\sqrt{2}x}+O(e^{-2\sqrt{2}x})
\end{align*}
 as $x\rightarrow\infty$.

By substituting the above asymptotic expansions into (\ref{eq:gamma_tau})
we obtain
\[
\frac{\partial}{\partial\gamma_{i}}\log\tau^{w}(x_{1},x_{2})=-\frac{\gamma_{i}}{4}\log x_{1}-\sum_{k=0}^{1}\frac{\gamma_{k}}{4}\left(\rho_{k}\right)_{\gamma_{i}}-\frac{\gamma_{i}}{4}+O(x_{1}^{\varepsilon_{i}}\log x_{1})+O(x_{2}^{\frac{3}{2}}e^{-2\sqrt{2}x_{2}})
\]
as $x_{1}\rightarrow0,\,x_{2}\rightarrow\infty$.

In our situation we have:

\[
\tau^{w}(1,x)=C_{0}x^{\frac{1}{8}(\gamma_{0}^{2}+\gamma_{1}^{2})}(1+O(x^{\varepsilon})),\quad x\rightarrow0,
\]

\[
\tau^{w}(1,x)=C_{\infty}e^{-x^{2}}(1+O(x^{1/2}e^{-2\sqrt{2}x})),\quad x\rightarrow\infty.
\]
 Then 
\[
\log\tau^{w}(x_{1},x_{2})=\log\frac{C_{\infty}}{C_{0}}-x_{2}^{2}-\frac{1}{8}(\gamma_{0}^{2}+\gamma_{1}^{2})\log x_{1}+O(x_{1}^{\varepsilon})+O(x_{2}^{1/2}e^{-2\sqrt{2}x_{2}}).
\]

Let 
\[
C:=\log\frac{C_{\infty}}{C_{0}}=\lim_{\substack{x_{1}\rightarrow0\\
x_{2}\rightarrow\infty
}
}\left(\log\tau^{w}(x_{1},x_{2})+x_{2}^{2}+\frac{\gamma_{0}^{2}+\gamma_{1}^{2}}{8}\log x_{1}\right).
\]

Then we obtain 
\[
\frac{\partial C}{\partial\gamma_{i}}=\lim_{\substack{x_{1}\rightarrow0\\
x_{2}\rightarrow\infty
}
}\left(\frac{\partial}{\partial\gamma_{i}}\left(\log\tau^{w}(x_{1},x_{2})+x_{2}^{2}+\frac{\gamma_{0}^{2}+\gamma_{1}^{2}}{8}\log x_{1}\right)\right)=-\frac{\gamma_{i}}{4}-\sum_{k=0}^{1}\frac{\gamma_{k}}{4}\left(\rho_{k}\right)_{\gamma_{i}}
\]
 that is, 

\begin{equation}
C=-\sum_{i=0}^{1}\frac{\gamma_{i}^{2}}{8}-\frac{1}{4}\sum_{k=0}^{1}\gamma_{k}\rho_{k}+\frac{1}{4}\int\sum_{k=0}^{1}\rho_{k}d\gamma_{k}.\label{eq:C}
\end{equation}

Note that $\frac{\partial K}{\partial m_{i}}=-2\frac{\partial K}{\partial\gamma_{i}}$,
where $K$ is the function defined in the proof of theorem \ref{thm:gen_fcn},
and
\[
\int\sum_{k=0}^{1}\rho_{k}d\gamma_{k}=-(\log2)\sum_{k=0}^{1}\gamma_{k}^{2}-2K+\text{{\rm const.}}
\]

The constant above is independent of the $\gamma_{i}$'s. By substituting
$\gamma_{0}=\gamma_{1}=0$, which corresponds to the trivial solution
$w_{0}\equiv w_{1}\equiv0$, into (\ref{eq:C}), we obtain $C=-4\left(\psi^{(-2)}(1/4)+\psi^{(-2)}(2/4)+\psi^{(-2)}(3/4)\right)+\text{{\rm const.}}$
On the other hand, the tau function $\tau^{w}(x_{1},x_{2})$ corresponding
to the trivial solution is $\exp(x_{1}^{2}-x_{2}^{2})$, so $C=0$
in this case. In conclusion we have the following result:
\begin{thm}
\label{thm:const}
\[
C=-\frac{1}{8}\left(\gamma_{0}^{2}+\gamma_{1}^{2}\right)-\frac{1}{2}\left(\gamma_{0}\rho_{0}+\gamma_{1}\rho_{1}\right)-\frac{1}{2}F+4\left(\psi^{(-2)}(1/4)+\psi^{(-2)}(2/4)+\psi^{(-2)}(3/4)\right).
\]
\end{thm}

The function $F$ in the theorem, which is the generating function,
is given in Definition \ref{def:generating function}.

\appendix

\section{Monodromy data\label{sec:Monodromy-data}}

At $\zeta=0$ we have a formal solution $\Psi_{f}^{(0)}=e^{-w}\Omega(I+\sum_{i\geq1}\Psi_{i}^{(0)}\zeta^{i})e^{\frac{1}{\zeta}d_{n+1}}$
of (\ref{eq:aux}), where 
\[
d_{n+1}=\begin{pmatrix}1\\
 & \omega\\
 &  & \ddots\\
 &  &  & \omega^{n}
\end{pmatrix},\;\Omega=\begin{pmatrix}1 & 1 & 1 & \cdots & 1\\
1 & \omega & \omega^{2} & \cdots & \omega^{n}\\
1 & \omega^{2} & \omega^{4} & \cdots & \omega^{2n}\\
\vdots & \vdots & \vdots & \ddots & \vdots\\
1 & \omega^{n} & \omega^{2n} & \cdots & \omega^{n^{2}}
\end{pmatrix}.
\]

We define the sector

\[
\Omega_{1}^{(0)}:=\begin{cases}
(-(\frac{1}{n+1}+\frac{1}{2})\pi,\frac{\pi}{2}) & (n+1\in2\mathbb{Z})\\
(-(\frac{1}{2(n+1)}+\frac{1}{2})\pi,(\frac{1}{2(n+1)}+\frac{1}{2})\pi) & (n+1\in2\mathbb{Z}+1)
\end{cases},
\]
where we use the notation $(a,b):=\{\zeta\in\mathbb{C}^{*}|a<\arg\zeta<b\}$.

We let $\Omega_{k+\frac{1}{n+1}}^{(0)}=e^{-\frac{\pi}{n+1}\sqrt{-1}}\Omega_{k}^{(0)}\:(k\in\frac{1}{4}\mathbb{Z})$
in the universal covering $\tilde{\mathbb{C}^{*}}$.

Let $\Psi_{k}^{(0)}$ be the fundamental solution such that $\Psi_{k}^{(0)}\sim\Psi_{f}^{(0)}$
on $\Omega_{k}^{(0)}$.

Similarly, at $\zeta=\infty$, we have the formal solution $\Psi_{f}^{(\infty)}=e^{w}\Omega^{-1}(I+\sum_{i\geq1}\Psi_{i}^{(\infty)}\zeta^{-i})e^{x^{2}\zeta d_{n+1}}$
and the sectors 
\[
\Omega_{1}^{(\infty)}:=\begin{cases}
(-\frac{\pi}{2},(\frac{1}{n+1}+\frac{1}{2})\pi) & (n+1\in2\mathbb{Z})\\
(-(\frac{1}{2(n+1)}+\frac{1}{2})\pi,(\frac{1}{2(n+1)}+\frac{1}{2})\pi) & (n+1\in2\mathbb{Z}+1)
\end{cases}
\]

\[
\Omega_{k+\frac{1}{n+1}}^{(\infty)}:=e^{\frac{\pi}{n+1}\sqrt{-1}}\Omega_{k}^{(\infty)}.
\]

Let $\Psi_{k}^{(\infty)}$ be the fundamental solution such that $\Psi_{k}^{(\infty)}\sim\Psi_{f}^{(\infty)}$
on $\Omega_{k}^{(\infty)}$.

We define the Stokes matrices $S_{k}^{(0)}$, $S_{k}^{(\infty)}$
by $\Psi_{k+1}^{(0)}=\Psi_{k}^{(0)}S_{k}^{(0)}$, $\Psi_{k+1}^{(\infty)}=\Psi_{k}^{(\infty)}S_{k}^{(\infty)}$.

We define the Stokes factors $Q_{k}^{(0)}$, $Q_{k}^{(\infty)}$ by
$\Psi_{k+\frac{1}{n+1}}^{(0)}=\Psi_{k}^{(0)}Q_{k}^{(0)}$, $\Psi_{k+\frac{1}{n+1}}^{(\infty)}=\Psi_{k}^{(\infty)}Q_{k}^{(\infty)}$.

Let $M:=Q_{1}^{(0)}Q_{1+\frac{1}{n+1}}^{(0)}\Pi$ where 
\[
\Pi=\begin{pmatrix}0 & 1\\
 & 0 & \ddots\\
 &  & \ddots & 1\\
1 &  &  & 0
\end{pmatrix}.
\]
Let $m_{i}$ be the eigenvalues of $M$. It is proved in \cite{II,III}
that the $m_{i}$ determine all $Q_{k}^{(0)}$.

We define the connection matrices $E_{k}$ by $\Psi_{k}^{(\infty)}=\Psi_{k}^{(0)}E_{k}$.

Let $E_{1}^{\text{global}}$ be $\frac{1}{n+1}AQ_{\frac{n}{n+1}}^{(\infty)}$
where 
\[
A=\begin{pmatrix}1\\
 &  &  & 1\\
 &  & \iddots\\
 & 1
\end{pmatrix}.
\]
 It is known that $E_{1}=E_{1}^{\text{global}}$ for global solutions
$w$. Let $e_{i}^{\mathbb{R}}$ be the eigenvalues of $E:=E_{1}\left(E_{1}^{\text{global}}\right)^{-1}$.
We have $e_{i}^{\mathbb{R}}e_{n-i}^{\mathbb{R}}=1$. It is proved
in \cite{III} that the $e_{i}^{\mathbb{R}}$ determine $E_{1}$.

\end{document}